# STOCHASTIC PROCESSES IN NETWORKS OF QUEUES WITH EXPONENTIAL SERVICE TIMES AND ONLY ONE CLASS OF CUSTOMERS

**M. A. M. Ferreira**

Instituto Universitário de Lisboa (ISCTE – IUL), BRU - IUL, Lisboa, Portugal

**Abstract**

Two networks of queues models, presented initially by Jackson, in the open case, and Gordon and Newell, in the closed case, stochastic processes are presented and studied in some of their details and problems. The service times are exponentially distributed and there is only one class of customers.

## 1 Introduction

In this paper, two networks of queues models, very simple, are studied. They were initially presented by Jackson, in the case of open networks [3,4] and by Gordon and Newell, in the case of closed networks[1] [12] . Posterior works, important on this subject are those of Baskett, Chandy, Muntz and Palacios [1] and of Kelly [2].

The simplicity of these models allows a detailed study of its behavior. So will be considered networks of queues with $J < \infty$ (finite number of nodes) and $L_j = \infty, j = 1,2,...,J$ (the whole nodes have infinite capacity). The service times at node $j$ are independent and exponentially distributed with parameter $\mu_j(n_j) > 0, n_j > 0$. This parameter may depend on the node and, for a fixed node, on the customers there present. It is supposed that $\mu_j(0) = 0, j = 1,2,...,J$ . Finally, the whole customers are of the same type.

## 2 Open Networks

For the open networks is still supposed that:

---

[1] A network is open if any customer can enter it or leave it. A network is closed if it has a fixed member of customers moving from node to node, with neither exterior arrivals nor departures. The networks opened for some customers and closed for others are called mixed, see for instance [5] .

.

-The exogenous arrivals process at node $j$ is a Poisson process with parameter $v_j$ that depends only on $j, j = 1,2,...,J$,

-The total exogenous arrivals process is a collection of independent Poisson processes or, in an equivalent way, there is only one arrivals process, choosing each customer, with probability $v_j / \sum_{k=0}^{J} v_k$, the node to enter in the network.

When abandoning the node, a customer is directed to another node or to the outside of the network in accordance with a commutation process.

**Definition 2.1**

Call $Y_m$ a random variable assuming values at $U'^2$, such that $Y_m = j$ if the node $j$ is the $m^{th}$ being visited by a given costumer, $\Im = \{Y_m : m = 1,2,...\}$ is the commutation process for that customer. ∎

Suppose that

- $\Im$ is a homogenous finite Markov chain,

- The path of every customer in the network is commanded by $\Im$.

The $\Im$ transition probabilities, called the commutation probabilities, are

$$P[Y_m = k | Y_{m-1} = j] = p(j,k), \quad k = 1,2,...,J, \quad j = 1,2,...,J$$
$$p(k) = P[Y_m = \Delta | Y_{m-1} = k] = 1 - \sum_{j=1}^{J} p(k,j), k = 1,2,...,J \quad (2.1).$$
$$P[Y_m = \Delta | Y_{m-1} = \Delta] = 1$$

The commutation probabilities depend on the node now occupied by the costumer and allow that, in the instant at which the commutation decision is taken, it is independent from everything happening in the network. The node $\Delta$ is an absorbent one.

It is also supposed that

- $\Im$ has no absorbent subsets in $U$,

-Every state of $U$ may be visited after any other,

-Every costumer may abandon the network,

-A costumer cannot go back immediately to the node from where it just left: $p(j,j) = 0, \forall j \in U$.

A consequence of this model, see for instance [10,11], is

---

[2] $U' = \{\Delta, 1,2,...,J\}$, where $\Delta$ is a node representing the exterior of the network, and $U = \{1,2,...,J\}$.

**Lemma 2.1**

The equations

$$\alpha_j = v_j + \sum_{k=1}^{J} \alpha_k + p(k,j), j = 1,2,\ldots,J \quad (2.2)$$

called network traffic equations, have only one solution $\alpha = (\alpha_1, \alpha_2, \ldots, \alpha_j), \alpha_j > 0, \forall j \in U$. ∎

Be $\mathcal{N} = \{N_1(t), N_2(t), \ldots, N_J(t), t \geq 0\}$ where $N_j(t)$ is the queue length in node $j, j = 1,2,\ldots,J$. If $\mathcal{N} = n = (n_1, n_2, \ldots, n_J)$ at instant $t$ it is said that the process $\mathcal{N}$ is in state $n$ at instant $t$.

It is supposed that $\mathcal{N}$ is stationary. Its states space is $\mathbb{N}_0^J$.

The stochastic process $\mathcal{N}$ study is easier if the following operators, defined in $\mathbb{N}_0^J$, are introduced, see [2, pg. 40 and pg. 48] and [10]:

$$T_{jk}n = (n_1, n_2, \ldots, n_j - 1, \ldots, n_k + 1, \ldots, n_J), \quad j < k, n_j > 0,$$
$$T_{jk}n = (n_1, n_2, \ldots, n_k + 1, \ldots, n_j - 1, \ldots, n_J), \quad j > k, n_j > 0$$
$$T_{j\cdot}n = (n_1, n_2, \ldots, n_j - 1, \ldots, n_J), \quad n_j > 0$$
$$T_{\cdot k}n = (n_1, n_2, \ldots, n_k + 1, \ldots, n_J)$$
$$(2.3).$$

Then, see [7],

**Lemma 2.2**

$\mathcal{N}$ is a homogenous continuous time Markov process, which transition rates are

$$q(n, T_{jk}n) = p(j,k)\mu_j(n_j),$$
$$q(n, T_{j\cdot}n) = p(j)\mu_j(n_j),$$
$$q(n, T_{\cdot k}n) = v_k,$$
$$q(n,n) = -\sum_{m \neq n} q(n,m), \quad (2.4). \blacksquare$$
$$q(n,m) = 0 \text{ in the remaining cases}$$

**Observation**:

-The only state changes that may occur in the network are the ones formalized by the operators $T_{jk}, T_{j\cdot}$ and $T_{\cdot k}$,

- The matrix, which elements are the $q(n,m)$ given by (2.4), is the $\mathcal{N}$ Markov process generator,

- The equation

$$\pi Q = 0 \qquad (2.5)$$

is the global equilibrium equation.

-If there is a solution of (2.5), $\pi$, such that $\pi > 0$ and $\sum_{n \in \mathbb{N}_0^J} \pi(n) = 1$, it is unique and it is called the $\mathcal{N}$ Markov process stationary distribution,

- The equation (2.5) in components form is

$$\pi(n)\left[\sum_j q(n, T_j.n) + \sum_j \sum_k q(n, T_{jk}n) + \sum_k q(n, T_{\cdot k}n)\right] =$$
$$\sum_j \pi(T_j.n) q(T_j.n, n) + \sum_j \sum_k \pi(T_{jk}n)q(T_{jk}n, n) +$$
$$\sum_k \pi(T_{\cdot k}n)q(T_{\cdot k}n, n) \qquad (2.6).$$

**Theorem 2.1**

The $\mathcal{N}$ Markov process stationary distribution $\pi$ is

$$\pi(n) = \prod_{j=1}^{J} \pi_j(n_j) \text{ where}$$

$$\pi_j(n_j) = \frac{b_j \, \alpha_j^{n_j}}{\prod_{r=1}^{n_j} \mu_j(r)} \text{ and} \qquad (2.7).$$

$$b_j = \sum_{i=0}^{\infty} \frac{\alpha_j^i}{\sum_{r=1}^{i} \mu_j(r)}$$

**Demonstration**: Note that (2.6) is fulfilled since the equations

$$\pi(n)\left[q(n, T_j.n) + \sum_k q(n, T_{jk}n)\right] = \pi(T_j.n)q(T_j.n, n) +$$

$$\sum_k \pi(T_{jk}n)q(T_{jk}n, n), \qquad j = 1,2, \dots, J \qquad (2.8)$$

and

$$\pi(n) \sum_k q(n, T_{\cdot k}n) = \sum_k \pi(T_{\cdot k}n) \, q(T_{\cdot k}n, n) \qquad (2.9)$$

so are. Substituting directly it is found that (2.7) satisfies (2.8) and (2.9). ∎

Observation:

- $\pi(n)$ is the nodal probabilities $\pi_j(n_j)$ product. This product form solution states that the vector process $\mathcal{N}$ coordinates are mutually independent for each $t$.

But, evidently, it is not allowed to conclude that $\mathcal{N}$ in $t_1$ and $\mathcal{N}$ in $t_2$ are mutually independent, although their stationary probabilities are identical.

- Suppose that $\pi_j(n_j) = \begin{cases} \mu_j, n_j > 0 \\ 0, n_j = 0 \end{cases}, j = 1,2,\dots,J$. So

$$\pi_j(n_j) = (1-\rho_j)\rho_j{}^{n_j}, n_j = 0,1,\dots \text{ and } \rho_j = \frac{\alpha_j}{\mu_j} < 1 \quad (2.10).$$

It seems, attending to the signification of the $\alpha_j$, that these are queue lengths generated by $M|M|1$ systems. But it is not completely like this, see [6,7].

Call $\mathcal{N}^r$ the $\mathcal{N}$ reverse process. If $q(n,n')$ are the $\mathcal{N}$ transition rates and $\pi(n)$ its stationary probabilities, the $\mathcal{N}^r$ transition rates are given by, see [2, pg. 28 and pg. 51],

$$q'(n',n) = \frac{\pi(n)q(n,n')}{\pi(n')} \quad (2.11).$$

Consequently,

### Lemma 2.3

The $\mathcal{N}^r$ transition rates are

$$q'(n, T_{jk}n) = \frac{\alpha_k \, p(k,j)\mu_j(n_j)}{\alpha_j},$$
$$q'(n, T_{j\cdot}n) = \frac{\nu_j \mu_j(n_j)}{\alpha_j},$$
$$q'(n, T_{\cdot k}n) = \alpha_k p(k),$$
$$q'(n,n) = -\sum_{n'\neq n} q(n,n'), \quad (2.12). \blacksquare$$
$$q'(n,n') = 0 \text{ in the remaining cases}$$

Then it is easy to conclude that

**Theorem 2.2**

If $\mathcal{N}$ is a stationary Markov process the same happens with $\mathcal{N}^r$. ∎

In the process $\mathcal{N}^r$, the arrivals from the outside of the network at the queue $k$ are a Poisson process at rate $\alpha_k p(k), k = 1, 2, \dots, J$. But these arrivals correspond precisely to departures from the system in the process $\mathcal{N}$. So

**Theorem 2.3**[3]

The departure processes from the various nodes of the network are independent Poisson processes, being the departure process from the node $k$ rate $\alpha_k p(k), k = 1, 2, \ldots, J$. If $t_0$ is a fixed instant, the departure process till $t_0$ and the process $\mathcal{N}$ state in $t_0$ are independent. ∎

# 3 Closed Networks

For closed networks

- The number of customers in the network is fixed and designated $N, N \geq 1$,
- $\mathcal{N}$ is still a Markov process with rates

$$q(n, T_{jk}n) = p(j,k)\mu_j(n_j) \quad (3.1)$$

and $T_{j\cdot}$ and $T_{\cdot k}$ have no meaning now. It is supposed that $\mathcal{N}$ is irreducible. Its states space is

$$\{(n_1, n_2, \ldots, n_J) \in \mathbb{N}_0^J : n_1 + n_2 + \cdots + n_J = N\} \quad (3.2).$$

**Theorem 3.1**

The $\mathcal{N}$ stationary distribution, $\pi$ is

$$\pi(n) = B_N \frac{\prod_{j=1}^{J} \alpha_j^{n_j}}{\prod_{r=1}^{n_j} \mu_j(r)} \quad (3.3)$$

where $B_N$ is the normalizing constant and $\alpha_j$ the total arrivals rate at node $j, j = 1, 2, \ldots J$, that must satisfy the traffic equations:

$$\alpha_j \sum_k p(j,k) = \sum_k \alpha_k \, p(k,j), j = 1, 2, \ldots, J \quad (3.4).$$

**Demonstration**: It is demonstrated such as Theorem (2.1) being now the global equilibrium equation

$$\pi(n) \sum_j \sum_k q(n, T_{jk}n) = \sum_j \sum_k \pi(T_{jk}n) q(T_{jk}n, n) \quad (3.5). \blacksquare$$

---

[3] The reverse process allows often the study, in a simple way, of the original process properties. Theorem 2.3 illustrates very well this idea, exhaustively used in [2].

Observation:

- The queue lengths are not independent since $B_N$ cannot be a product of $b_j$'s because $n_1 + n_2 + \cdots + n_J = N$. In spite of this (3.3) is called also a product form solution.

It is usually very difficult and fastidious to compute the $B_N$ value for closed queuing networks. This difficulty increases with the network number of nodes. Two results of Harrison, see [9], that allow the use of more efficient algorithms, are then presented.

If $\mu_j(n_j) = \mu_j, j = 1,2,\ldots,J$ (3.3) may be written as

$$\pi(n) = B_N \prod_{j=1}^{J} \rho_j^{n_j} \qquad (3.6)$$

being $\rho_j$ given as in (2.10).

**Theorem 3.2**

If the $\rho_j, j = 1,2,\ldots,J$ values are all different

$$B_N = \frac{\sum_{j=1}^{J} \rho_j^{N+J-1}}{\prod_{i \neq j}(\rho_j - \rho_i)} \qquad (3.7).$$

**Demonstration**: It is made by induction in $J$, after checking that it is true for $J = 1$. ∎

**Theorem 3.3**

If the $\rho_j$ are not all different (degenerate case)

$$B_N = \left\{ \prod_{j=m+1}^{J} (1 + \rho_j D_j)^{d_j} \right\} b_m \qquad (3.8)$$

where

-$d_j$, associated to $\rho_j$, is the cardinal of the set

$$\{i: 1 \leq i \neq j \leq J, \rho_i = \rho_j\} \qquad (3.9),$$

-$b_m$ is the normalizing constant for the sub network constituted by the $j$ nodes in the set $\{m_1, m_2, \ldots, m_j\}$ with the same vector $\{\rho_{m_1}, \rho_{m_2}, \ldots, \rho_{m_j}\}$, where for the whole $1 \leq k \leq j, m_k = i$ for some $1 \leq i \leq J$,

-$D_j$ means derivation to $\rho_j$.

**Demonstration**: It is driven like the one of Theorem 3.2, with the necessary adaptations to the degeneration.∎

## 4 Conclusions

The networks of queues models considered in this work are the most popular in the practical applications.

The analytical forms of the population process solutions induce some confusion in what concerns statistical independence, with influence in the so-called product form solutions, and in the characterization of the nodes behavior, as if it was the same as for the isolated nodes. This was clarified in the above text and this only may be achieve studying the stochastic processes of the networks. The clarification of this situations is important because if it is not done the networks parameters estimation may be prejudiced.

The problem of the normalizing constant computation for closed networks was also presented and solutions for it were given.

Often, in queue systems study, it is very difficult to obtain analytical friendly expressions for the various parameters and quantities of interest. So it is usual to make use of numerical methods and simulation[4]. In the simulation case it is evident the importance of knowing very deeply the stochastic processes involved. If this does not happen it is impossible to make a correct simulation of the system. This is one more and important reason to study carefully these stochastic processes.

---

[4] On simulation of queue systems see, for instance, [8] .